\documentclass[11pt]{article}
\usepackage{amscd}      
\usepackage{amssymb}
\usepackage[intlimits]{amsmath}

\usepackage{theorem}

\newtheorem{prop}{Proposition}[section]
\newtheorem{lem}[prop]{Lemma}

\newtheorem{them}[prop]{Theorem}

\theorembodyfont{\upshape}

\newtheorem{defn}[prop]{Definition}

\newtheorem{numex}[prop]{Example}

\newtheorem{rmk}{Remark}

\newenvironment{pf}{\begin{trivlist}\item[]{\sc Proof.}}%
           {\nolinebreak $\Box$ \end{trivlist}}

\newcommand{\noprint}[1]{}

\renewcommand{\tilde}{\widetilde}

\newcommand{\alg}{\mbox{\tiny alg}}

\newcommand{\rr}{{\mathbb R}}

\def\lcf{\lbrack\! \lbrack}
\def\rcf{\rbrack\! \rbrack}

\def\alg{(A,\lcf \cdot , \cdot \rcf ,\rho )}
\def\gfJ{graph \, J}
\def\S{\mathbb{S}}

\newcommand{\gm}{\Gamma }

\newcommand{\backl}{\mathbin{\vrule width1.5ex height.4pt\vrule height1.5ex}}
\newcommand{\per}{\backl}
\newcommand{\be}{\begin{eqnarray*}}
\newcommand{\ee}{\end{eqnarray*}}

\newcommand{\half}{\frac{1}{2}}

\begin{document}
\title{{\bf Hamiltonian spaces for Manin pairs over manifolds}}
\author{David Iglesias Ponte$^1$, Ping
Xu$^2$
\\{\small\it $^1$Instituto de Ciencias Matem\'aticas (CSIC-UAM-UCM-UC3M),
Spain}
\\{\small\it $^2$Department of Mathematics, Pennsylvania State
University, USA}\\[5pt]
{\small\it e-mail: iglesias@imaff.cfmac.csic.es,
ping@math.psu.edu} }
\date{}

\sloppy
\maketitle

\begin{abstract}
We introduce the notion of  Hamiltonian spaces for Manin pairs over
manifolds, using the so-called generalized Dirac structures. As an example,
we describe  Hamiltonian spaces of a  quasi-Lie bialgebroid
using this general framework. We also discuss
reduction of Hamiltonian spaces of this general type.
\end{abstract}

\section{Introduction}

In recent years, several Lie-group valued moment map theories
haven been introduced, extending the usual  Hamiltonian
action of a Lie group $G$ on a symplectic manifold $M$ (in this case, the target of the
moment map is $\mathfrak g^*$, the dual space of the Lie algebra
$\mathfrak g$ of $G$), including McDuff's $\S^1$ valued moment
maps \cite{M}, Lu's Poisson Lie group valued moment
maps \cite{LuW:1990} and the Lie group valued moment maps of
Alekseev {\em et al} \cite{AMM}.

\medskip
In order to give a common framework to all of these, in \cite{AK}
Alekseev and Kosmann-Schwarzback introduced  a moment map theory
 based on Manin
pairs $(\mathfrak d, \mathfrak g)$ and Manin quasi-triples. In
such a  case, the moment map takes values in the homogeneous space $D/
G$, where $D$ and $G$
are  Lie groups integrating $\mathfrak d$ and
 $\mathfrak g$, respectively. More precisely, given
a Manin quasi-triple $(\mathfrak d, \mathfrak g,\mathfrak h)$,
$\mathfrak h$ being  an isotropic complement of
a maximal isotropic Lie subalgebra $\mathfrak g$ in
$\mathfrak d$,
 one constructs a Lie quasi-bialgebra $(\mathfrak g,
F, \Omega )$ which can be integrated to a quasi-Poisson Lie group.
A quasi-Poisson space is a $G$-manifold $X$ endowed with a
bivector $\Pi _X$ satisfying certain compatibility conditions. It
is quasi-hamiltonian if there exists  a moment map $J:X\to D/G$,
which induces the $G$-action on $X$ in a certain sense.
 In  particular, if  $(\mathfrak d, \mathfrak g,\mathfrak h)$ is a
 Manin triple, that is,
$(\mathfrak g,\mathfrak h)$ defines a Drinfeld's Lie bialgebra,
we recover Lu's moment map for Poisson group actions.


\medskip
The notion of Courant algebroids is a natural  generalization of double of
Lie bialgebras  \cite{LWX0}.
 Roughly speaking, a Courant algebroid is a vector bundle
 $E\to S$ endowed with a non-skew symmetric
bracket on $\gm (E)$, a bundle map $\tilde{\rho}:E\to TS$ and a nondegenerate
bilinear form on $E$, satisfying certain axioms similar to those of a Lie
algebroid.
 One can also define the  notion of
Dirac structures.
  When the base space $S$ is a single
point, a Courant algebroid with  a Dirac structure
reduces to  a Manin pair $(\mathfrak d,
\mathfrak g)$.
 Courant algebroids and Dirac structures have been used
in the   study of  Poisson homogeneous
spaces for Poisson groupoids and Poisson reduction (for more
details, see \cite{LWX}).

\medskip
Closely related to Courant algebroids is the notion of quasi-Lie
bialgebroids, first introduced by Roytenberg \cite{Roy}, which is
a natural generalization of Lie bialgebroids
\cite{MackenzieX:1994}. It also generalizes Drinfeld's quasi-Lie
bialgebras \cite{Drinfeld}, the classical limit of quasi-Hopf
algebras. More precisely, given a Courant algebroid $E\to S$ and a
Dirac structure $A$ on it, if we choose an isotropic complement
$B$ of $A$ in $E$,
we obtain a quasi-Lie bialgebroid. Conversely, from a quasi-Lie
bialgebroid structure on a Lie algebroid $A\to S$,  one
can construct a Courant algebroid structure on $A\oplus A^*\to S$.
The correspondence is analog to the
that between Manin quasi-triples and Lie
quasi-bialgebras. In addition, if the isotropic complement
$B$  is integrable, that is, it is also a Dirac structure, we
have a Manin triple of a Lie bialgebroid \cite{LWX0}.

\medskip
The main purpose of this paper is to introduce the notion of a
Hamiltonian space for a Manin pair over a manifold, that is, a
Courant algebroid ${\mathfrak D}\to S$ along with a Dirac
structure $A\to S$ on it. In order to do so, we first recall the
notion of a generalized Dirac structure for a Courant algebroid
(see \cite{AX}). This is a Lagrangian subbundle $F\to Q$ of the Courant
algebroid where the integrability conditions have been weakened,
since the base space $Q$ is only a  submanifold of $S$. More
precisely, a Hamiltonian space for a Manin pair $(\mathfrak D ,A)$
is a manifold $X$ together with a smooth map $J:X\to S$ such that
there exists a generalized Dirac structure $F\to graph\, J$ in the
product Courant algebroid ${\mathfrak D}\times (TX\oplus T^*X)\to
S\times X$ satisfying some extra conditions. Our work
was inspired by \cite{BCS}, where  quasi-Poisson spaces for
a Manin pair $(\mathfrak d,\mathfrak g)$
are interpreted as Dirac structures on
$\mathfrak d \times (TX\oplus T^*X) \to X$.

\medskip
The paper is organized as follows. In Section 2 we review the
definition of a Courant algebroid and its relation with quasi-Lie
bialgebroids. We illustrate the notion with some examples which
will be useful along the paper, in particular, the quasi-Lie
bialgebroid associated with a Manin-quasi triple. Finally, we
recall the notion of a generalized Dirac structure for a Courant
algebroid, the notion in which our  main definition  is based on.
In Section 3, we introduce the notion of a Hamiltonian
space for a Manin pair $(\mathfrak D, A)$ over a manifold $S$ with
moment map $J:X\to S$. From this concept it is possible to deduce
an action of the Lie algebroid $A$ on $J:X\to S$ and a Dirac
structure $F(A)$ on $X$. This Dirac structure will allow us to
apply the standard reduction procedure \cite{Co}
 and construct a Poisson structure on  its
 quotient space $X/_{\widetilde{ }}\,$. In Section 4, we
study the particular situation when we choose an isotropic
complement $B$ of $A$ in $\mathfrak D$. More precisely, we prove
that Hamiltonian spaces for quasi-Lie bialgebroids are in
one-to-one correspondence with Hamiltonian spaces for the
corresponding Manin pair. Several examples are discussed: Poisson
actions of Lie-bialgebroids on Poisson manifolds and
quasi-Hamiltonian spaces for Manin quasi-triples.

\medskip
Recently, Dirac realizations have been used to characterize
quasi-hamiltonian spaces \cite{BC,BC2}.
The connection between \cite{BC,BC2}
and our theory is explored in
  a subsequent paper \cite{BI}.

\section{Preliminaries}
We start recalling several notions which will be used along the
paper.
\subsection{Lie algebroids}
A \textbf{Lie algebroid} $A$ over a manifold $S$ is a vector bundle
$A$ over $S$ together with a Lie bracket $\lcf\cdot ,\cdot \rcf$ on
the space of sections $\Gamma (A)$
and a bundle map $\rho \colon A \to TS$, called the \textbf{anchor map},
such that $\rho$ induces a Lie algebra morphism
$\Gamma (A) \to \mathfrak X (S)$ and
 the following compatibility condition holds:
\[
\lcf X ,fY \rcf=f\lcf X,Y\rcf +(\rho  (X )(f))Y.
\]
for all $f\in C^\infty (S)$ and for all $X ,Y \in \Gamma (A)$. The
triple $\alg$ is called a Lie algebroid over $S$ (see \cite{Mk}).

Let $\alg$ be a Lie algebroid over a manifold $S$ and $J :X\to S$ be
a smooth map. An \textbf{action} of $A$ on $J:X\to S$ is a
$\rr$-linear map
\[
\Gamma (A)\to \mathfrak X (X),\quad X\in \Gamma (A)\mapsto \hat{X}
\in \mathfrak X (X),
\]
such that:
\[
\widehat{fX} =(J^*f)\hat{X},\quad \widehat{\lcf X, Y\rcf} =[
\hat{X} ,\hat{Y} ],\quad J_\ast(\hat{X} (x))=\rho (X (J(x))),
\]
for $f\in C^\infty (S)$, $X,Y\in \Gamma (A)$ and $x\in X$.
This is simply an infinitesimal version of Lie groupoid action.

\subsection{Courant algebroids}

\begin{defn}
A \textbf{Courant algebroid} over a manifold $S$ is a vector
bundle $E\to S$ equipped with a nondegenerate symmetric bilinear
form $( \cdot | \cdot )$ on the bundle, a bundle map $\tilde\rho
:E\to TS$ and a bilinear bracket $\circ$ on $\Gamma
(E)$, called {\em Dorfman bracket},
 such that for all $e, e_1,e_2,e_3\in \Gamma (E)$, $f\in
C^\infty (S)$ the following  axioms are  satisfied:
\begin{itemize}
\item[\it i)] $ e_1\circ (e_2 \circ e_3) =( e_1 \circ e_2) \circ e_3
\} + e_2 \circ  (e_1 \circ e_3)$;
\item[\it ii)] $ e \circ e ={\cal D}( e\, |\, e )$;
\item[\it iii)] ${\cal L}_{\tilde\rho (e)}( e_1\, |\,e_2) =
( e\circ e_1 |\,e_2)+( e_1\, |\, e \circ e_2 )$;
\item[\it iv)] $\tilde\rho ( e_1 \circ e_2 )=[ \tilde\rho (e_1), \tilde\rho (e_2)]$;
\item[\it v)] $ e_1  \circ (f\, e_2) =f( e_1 \circ  e_2) +
 ({\cal L}_{\tilde\rho (e_1)}f)e_2$,
\end{itemize}
where ${\cal D}:C^\infty (S)\to \Gamma (E)$ is defined by $( {\cal
D} f\, |\, e)={\cal L}_{\tilde\rho (e)}f$.
Note that the Dorfman bracket $\circ$ is not skew-symmetric in general.
 Its skew-symmetrization is called {\em Courant bracket}
and satisfies a set of more complicated equations \cite{LWX}.

A subbundle $L\subset E$ is called a \textbf{Dirac structure} (or
a \textbf{Dirac subbundle}) if it is maximal isotropic with
respect to $( \cdot |\cdot )$ and if $\Gamma (L)$ is closed
 $\{ \cdot ,\cdot \}$.
\end{defn}

The following is well known.

\begin{prop}
If $L\subset E$ is a Dirac structure for the Courant algebroid $E\to
S$ then $L\to S$ is a Lie algebroid, where the bracket $\lcf \cdot
,\cdot \rcf$ is just the restriction of $\{ \cdot ,\cdot \}$ to
$\Gamma (L)$ and the anchor $\rho :L\to TS$ is the restriction of
$\tilde\rho: E\to TS$ to $L$.
\end{prop}
\begin{numex}
A \textbf{Manin pair} $({\mathfrak d}, {\mathfrak g})$ is an even
dimensional Lie algebra ${\mathfrak d}$ with an invariant,
nondegenerate symmetric bilinear form (which is a Courant algebroid
over a point) and $\mathfrak g$ is a maximal isotropic subalgebra
(Dirac structure) of $\mathfrak d$
\end{numex}
\begin{numex}
Let $M$ be an arbitrary manifold. On the vector bundle $E_M
=TM\oplus T^*M$ over $M$, there exists a natural symmetric pairing
$( \cdot |\cdot ) $ given by
\begin{equation}\label{pairings}
\begin{array}{c}
( (X_1,\alpha _1)\, |\,(X_2,\alpha _2)) = \frac{1}{2}( \alpha _1(X_2)+\alpha
_2(X_1)),
\end{array}
\end{equation}
for $(X_i,\alpha _i)\in \Gamma (TM\oplus T^*M)$, $i\in \{ 1,2\}$. On
the other hand, the \textbf{Courant bracket} $[\cdot ,\cdot ]_c$ on
$\Gamma (TM\oplus T^*M)$ is defined by
\begin{equation}\label{bracket}
[ (X_1,\alpha _1),(X_2,\alpha _2)]_c= \big ( [X_1,X_2] , {\cal
L}_{X_1}\alpha _2-i_{X_2}d\alpha _1  \big ),
\end{equation}
where $[\cdot ,\cdot ]$ is the usual Lie bracket of vector fields
and ${\cal L}$ is the Lie derivative operator on $M$. We have that
$(E_M, ( \cdot |\cdot ) ,[\cdot ,\cdot ]_c, \tilde{\rho})$ is a
Courant algebroid, where $\tilde{\rho}:E_M\to TM$ is the first
projection. This was first introduced in \cite{Co}. An important
feature of Dirac structure  lies in the study of reduction
\cite{Co}, which we recall briefly below.

Given a Dirac structure $L$ on a manifold $M$, a function $f\in
C^\infty (M)$ is said to be \textbf{admissible} if there exists
$X_f\in \mathfrak X(M)$ such that $(X_f,df)\in \Gamma (L)$. The
space of admissible functions, denoted by $C^\infty _L(M)$, is
endowed with a Poisson algebra structure as follows. If $f,g\in
C^\infty _L(M)$ then
\[
\{ f,g \}_L = X_f(g).
\]
A Dirac structure $L$ on a manifold $M$ induces a singular
 distribution  $L\cap TM$, which  is called
the \textbf{characteristic distribution} of  $L$.
In addition, if  $M / (L\cap TM)$ is a smooth
manifold, it is easy to see  that $C^\infty (M/ (L\cap TM))$ is
isomorphic $C^\infty _L(M)$. Thus $M/(L\cap TM)$ is a Poisson manifold.
\end{numex}

\subsection{Quasi-Lie bialgebroids and Courant algebroids}
Quasi-Lie bialgebroids, which were first introduced in \cite{Roy},
are closely related with Courant algebroids and Dirac structures.
Here, we give an alternative interpretation using 2-differentials
\cite{ILX}.

\begin{defn}\label{quasi-Lie-bialgebroid}
A \textbf{quasi-Lie bialgebroid} consists of
a triple   $(A, \delta ,\Omega )$, where $A$ is a Lie
algebroid, $\delta$ is a 2-differential, and
$\Omega \in \gm (\wedge^3 A )$  such that
\[
\delta ^2 =\lcf \Omega , \cdot \rcf ,\qquad \delta \Omega =0.
\]
Here, by a 2-differential on a Lie algebroid $A$, we mean
 a linear operator $\delta :\gm (\wedge ^\bullet A)\to \gm
(\wedge ^{\bullet +1}A)$ satisfying
\begin{itemize}
\item[{\it i)}] $\delta (P\wedge Q)=(\delta P)\wedge Q+(-1)^{p}P\wedge \delta
Q$,
\item[{\it ii)}] $\delta \lcf P, Q\rcf =\lcf \delta P, Q\rcf +(-1)^{(p+1)}
\lcf P, \delta Q\rcf$,
\end{itemize}
for all $P\in \gm (\wedge ^pA)$ and $Q\in\gm (\wedge ^q A)$.
\end{defn}

\begin{numex}
If $\Omega =0$,  then $\delta$ is a 2-differential of square $0$.
Thus, one can define a Lie algebroid structure on $A^*$
 by
\[
\begin{array}{lll}
\rho _\ast (\xi _1)(f)&=&\langle \delta
f,\xi  _1 \rangle , \ \ \mbox{ and}\\
[10pt] \langle \lcf \xi _1 , \xi _2 \rcf_\ast ,X\rangle &=&
 \rho _\ast (\xi _1 )
(\xi _2(X)) - \rho _\ast (\xi _2 ) (\xi _1(X))  - \langle \delta X,
\xi _1\wedge \xi _2\rangle ,
\end{array}
\]
for $\xi _1 ,\xi _2 \in \gm (A^\ast )$, $X\in \gm (A)$ and $f\in
C^\infty (S)$. Using \textit{ii)} in Definition
\ref{quasi-Lie-bialgebroid}, we conclude that the pair $(A,A^*)$ is
a \textbf{Lie bialgebroid} (see \cite{MackenzieX:1994}).
\end{numex}

\begin{numex} \cite{AK}\label{Manin-pair}
Let $(\mathfrak d, \mathfrak g)$ be a \textbf{Manin pair} and
$(D,G)$ the corresponding group pair, i.e. $D$ and $G$ are
connected and simply connected Lie groups with Lie algebras
${\mathfrak d}$ and $\mathfrak g$ respectively. Furthermore, the
action of the Lie group $D$ on itself by left multiplication induces
an action of $D$ on $S=D/G$,  and in particular a $G$-action on $S$,
which is called the \textbf{dressing action}. As in \cite{AK}, the
infinitesimal dressing action is denoted by $v\mapsto v_S$ for any
$v\in \mathfrak d$.

\medskip
If $\mathfrak h$ is an isotropic complement of $\mathfrak g$ in
$\mathfrak d$, by identifying $\mathfrak h$ with $\mathfrak g^\ast$,
we obtain a quasi-Lie bialgebra structure on $\mathfrak g$, with
the cobracket $F:\mathfrak g\to \wedge ^2 \mathfrak g$ and $\Omega \in
\wedge ^3 \mathfrak g$. If $\{ e_i \}$ is a basis of $\mathfrak g$
and $\{ \epsilon ^i\}$ the dual basis of $\mathfrak g^\ast \cong
\mathfrak h$, then we may use
 $F(e_i)=\half \sum _{j,k}F^{jk}_i e_j\wedge e_k$
and $\Omega =\frac{1}{6}\sum _{i,j,k} \Omega ^{ijk}e_i\wedge
e_j\wedge e_k$. Moreover, the bracket on $\mathfrak d \cong
\mathfrak g\oplus \mathfrak h$ can be written as
\begin{equation}\label{d-bracket}
[e_i,e_j]_{\mathfrak d}\kern-3pt=\kern-3pt\displaystyle\sum _{k=1}^n
c_{ij}^k e_k,\quad
 [e_i,\epsilon ^j]_{\mathfrak d}\kern-3pt=\kern-3pt\displaystyle\sum _{k=1}^n
-c^j_{ik}\epsilon ^k +F_i^{jk}e_k, \quad [\epsilon ^i,\epsilon
^j]_{\mathfrak d}\kern-3pt=\kern-3pt\displaystyle\sum _{k=1}^n
F_k^{ij}\epsilon ^k+ \Omega ^{ijk}e_k,
\end{equation}
where $c_{ij}^k$ are the structure constants of the Lie algebra
$\mathfrak g$ with respect to the basis $\{ e_i \}$. We call
$(\mathfrak d, \mathfrak g, \mathfrak h)$ a \textbf{Manin
quasi-triple}. Conversely, given a a quasi-Lie bialgebra $(\mathfrak
g, F,\Omega )$, then $(\mathfrak g\oplus \mathfrak g^*,\mathfrak g,
\mathfrak g^*)$ is a Manin-quasi triple, where the Lie algebra
structure on $\mathfrak g\oplus \mathfrak g^*$ is given by Eq.
(\ref{d-bracket}) and the pairing is the canonical one.

\medskip
Given $(\mathfrak d, \mathfrak g, \mathfrak h)$ a Manin
quasi-triple, let $\lambda :T^\ast _sS\to \mathfrak g$ be the dual
map of the infinitesimal dressing action $\mathfrak g^\ast \cong
\mathfrak h \to T_sS$. That is, $\langle \lambda (\theta _s),\eta
\rangle = \langle \theta _s,\eta _S(s) \rangle$, for $\theta _s \in
T^\ast _sS$ and $\eta \in \mathfrak h$. A direct consequence is
\begin{equation}\label{def-lambda}
\lambda (df)=\displaystyle \sum _{i=1}^n (\epsilon
^i)_S(f)e_i,\mbox{ for }f \in C^\infty (S).
\end{equation}

Assume that $(\mathfrak d, \mathfrak g, \mathfrak h)$ is a Manin
quasi-triple with associated quasi-Lie bialgebra $(\mathfrak
g,F,\Omega )$. On the transformation Lie algebroid $\mathfrak
g\times S\to S$ define an almost 2-differential
\begin{equation}\label{2-diff-trans}
\begin{array}{l}
\delta  (f)=\lambda (df),\mbox{ for }f \in
C^\infty (S),\\[10pt]
\displaystyle \delta  \xi = - F(\xi),\mbox{ for }\xi \in \mathfrak
g,
\end{array}
\end{equation}
where $\lambda$ is defined by Eq. (\ref{def-lambda}) and $\xi \in
\mathfrak g$ is considered as a constant section of the Lie
algebroid $\mathfrak g\times S\to S$, extending this operation to
arbitrary section using the derivation law. Then, $(\mathfrak
g\times S, \delta ,\Omega )$ is a quasi-Lie bialgebroid (for more
details, see \cite{ILX}).
\end{numex}
\bigskip Next, we will recall the relation between quasi-Lie bialgebroids
and Courant algebroids, which generalize the correspondence between
quasi-Lie bialgebras and Manin quasi-triples.

Let $(A,\delta ,\Omega )$ be a quasi-Lie bialgebroid. Then, using
$\delta$, one can construct a pair $(\lcf \cdot, \cdot \rcf _* ,\rho
_* )$ as follows. Let
\begin{equation}\label{bracket+anchor}
\begin{array}{lll}
\rho _\ast (\xi _1)(f)&=&\langle \delta
f,\xi  _1 \rangle , \ \ \mbox{ and}\\
[10pt] \langle \lcf \xi _1 , \xi _2 \rcf_\ast ,X\rangle &=&
 \rho _\ast (\xi _1 )
(\xi _2(X)) - \rho _\ast (\xi _2 ) (\xi _1(X))  - \langle \delta X,
\xi _1\wedge \xi _2\rangle ,
\end{array}
\end{equation}
for $\xi _1 ,\xi _2 \in \gm (A^\ast )$, $X\in \gm (A)$ and $f\in
C^\infty (S)$. Thus, we can obtain Courant algebroid $E=A\oplus A^*$
where the anchor $\tilde\rho$ and the bracket $\{ \cdot ,\cdot \}$
are given by
\begin{equation}\label{double+Courant}
\begin{array}{rcl}
\tilde\rho&=&\rho +\rho _*\\[7pt]
\{ (a_1,\xi _1),(a_2,\xi _2) \}&=& (\lcf a_1,a_2\rcf - i_{\xi
_2}\delta a_1 +{\cal L}_{\xi _1}a_2 +(\xi _1\wedge \xi _2)\per
\Omega , \\[3pt] & &\lcf \xi _1,\xi _2\rcf _* - i_{a _2}d_A \xi _1 +{\cal
L}_{a_1}\xi _2).
\end{array}
\end{equation}
Conversely, suppose that $E\to S$ is a Courant algebroid and that
$A\to S$ is a Dirac structure for $E$. If we choose $H$ a
complementary isotropic subbundle of $A$ (not necessarily
integrable) then one can identify $H$ with $A^*$ using $( \cdot
|\cdot )$ so that $E\cong A\oplus A^*$. In addition, the anchor
map and the bracket on $E\cong A\oplus A^*$ are defined by Eq.
(\ref{double+Courant}), where $\lcf \cdot ,\cdot \rcf _*$, $\rho
_*$ and $\Omega\in \gm (\wedge ^3 A)$ are characterized by
\[
\begin{array}{l}
\lcf \xi _1,\xi _2 \rcf _*=\mbox{pr}_{A^*}(\{ (0,\xi _1),(0,\xi _2)
\}
)\\[7pt]
\rho _* (\xi _1)=\tilde{\rho}((0,\xi _1))\\[10pt]
\Omega (\xi _1,\xi _2)=\mbox{pr}_A \{ (0,\xi _1),(0,\xi _2)\},
\end{array}
\]
for $\xi _1,\xi _2\in \gm (A^*)$. If we construct $\delta :\gm
(\wedge ^\bullet A)\to \gm (\wedge ^{\bullet +1}A)$ mimicking the
definition of the differential of a Lie algebroid, we see that
$(A,\delta ,\Omega )$ is a quasi-Lie algebroid (for more details,
see \cite{Roy}). If the isotropic complement $H\to S$ is also a
Dirac structure then, $\Omega =0$ and, as a consequence, $(A, A^*)$
is a Lie bialgebroid.
\subsection{Generalized Dirac structures for Courant algebroids}
In this Section, we will generalize the notion of a Dirac structure
$A$ for a Courant algebroid $E\to S$, allowing the base space $Q$ of
$A$ to be a submanifold of $S$.
\begin{defn}\label{gds}
{\rm \cite {AX} Given a Courant algebroid $E\to S$, a Lagrangian
subbundle $F\to Q$ over a submanifold $Q$ of $S$ is a
\textbf{generalized Dirac structure} if
\begin{itemize}
\item[{\it i)}] $F$ is compatible with the anchor, that is, $\tilde\rho (F)\subset
TQ$, where $\tilde\rho :E\to TS$ is the anchor of the Courant
algebroid;
\item[{\it ii)}] for any sections $e_1,e_2$ of $E$ such that
$e_1|_Q$,$e_2|_Q\in \gm (F)$, $\{ e_1,e_2\} |_Q\in \gm (F)$, where
$\{ \cdot ,\cdot \}$ is the Courant algebroid bracket.
\end{itemize}
}
\end{defn}

\begin{rmk}
In the above definition, the bracket on ii) can also be replaced
by the Dorfman bracket $\circ$. One can easily check that
they are equivalent. This is because
$e_1\circ e_2|_Q-\{ e_1,e_2\} |_Q=D(e_1, e_2)|_Q\in \Gamma (F)$
since $(D(e_1, e_2), e)=0$ for any $e\in \Gamma (F)$.
\end{rmk}

\begin{numex}
Let $(S,\pi)$ be a Poisson manifold and $TS\oplus _{\pi} T^*S$ be
the corresponding Courant algebroid coming from the Lie bialgebroid
structure $(TS,\delta _\pi,0)$, where $\delta _\pi:\gm (\wedge
^\bullet (TS))\to \gm (\wedge ^{\bullet +1}(TS))$ is given by
$\delta _\pi P=[\pi ,P]$. Then, it is not difficult to show that a
submanifold $Q$ of $S$ is a coisotropic submanifold if and only if
$TQ\oplus NQ\to Q$ is a generalized Dirac structure for $TS\oplus
_\pi T^*S$, $NQ$ being the normal bundle of $Q$.
\end{numex}

\begin{numex}\label{ej-gds}
Let $(A,A^*)$ and $(B,B^*)$ be Lie bialgebroids over bases $M$ and
$N$ respectively, and $E_1=A\oplus A^*$ and $E_2=B\oplus B^*$ be
their doubles. Assume that $\Phi :A\to B$ is a Lie bialgebroid
morphism over $\phi :M\to N$. Then
\[
F=\{ ((a,-\Phi ^*b_*),(\Phi a,b_*)) \,|\, \forall\, a\in A_x\mbox{
and } b_*\in B_{\phi (x)}^* \} \subset E_1\times E_2
\]
is a generalized Dirac structure.
\end{numex}
\section{Hamiltonian spaces for Manin pairs over manifolds}
In this Section, we will introduce the notion of Hamiltonian space
for a Manin pair over a manifold. Some of their properties will be
developed and a reduction process will be explained.
\subsection{Definition and main properties}
Let $(\mathfrak D, A)$ be a \textbf{Manin pair} over $S$, that is,
$\mathfrak D\to S$ be a Courant algebroid and $A\to S$ a Dirac
structure on it.
\begin{defn}\label{Hamiltonian-space}
A \textbf{Hamiltonian space} for the Manin pair $(\mathfrak D, A)$
(over the manifold $S$) is a manifold $X$ together with a smooth
map $J:X\to S$ such that there exists a generalized Dirac
structure $F\subset E=\mathfrak D\times (TX\oplus T^*X)$ over
$Q=\gfJ$ such that
\begin{itemize}
\item[{\it i)}] $F$ intersects $TX$ trivially;
\item[{\it ii)}] the intersection of $F$ with $ E\times (TX\oplus \{ 0\})$
projects onto $A$ under the natural map $E\to \mathfrak D$.
\end{itemize}
\end{defn}
\begin{lem}\label{vlemma}
Let $J:X\to S$ be a Hamiltonian space for the Manin pair
$(\mathfrak D, A)$. Then,
\begin{itemize}
\item[{\it i)}] For any $a\in A_{J(x)}$, there exists a unique $u\in
T_xX$ such that $(a,(u,0))\in F_{(J(x),x)}$.
\item[{\it ii)}] For any $\alpha \in T_x^*X$, there exists $d\in
{\mathfrak D}_{J(x)}$ and $u\in T_xX$ such that $(d,(u,\alpha ))\in
F_{(J(x),x)}$.
\end{itemize}
\end{lem}
\begin{pf}
Using {\it ii)} in Definition \ref{Hamiltonian-space}, we know that
for any $a\in A_{J(x)}$ there exists $u\in T_xX$ such that
$(a,(u,0))\in  F_{(J(x),x)}$. If $u'\in T_xX$ such that
$(a,(u',0))\in F_{(J(x),x)}$ then, we have that $(0, (u-u',0))\in
F_{(J(x),x)}$. Using {\it i)} in Definition \ref{Hamiltonian-space},
we conclude that $u=u'$. This proves \textit{i)}.

\medskip
Now, consider the projection $\varphi : F\to T^*X$. Using {\it ii)}
in Definition \ref{Hamiltonian-space}, one can deduce that
$\mbox{Ker }\varphi =F\cap ({\mathfrak D}\times (TX\oplus \{
0\}))\cong A$. Thus, $\varphi$ is surjective, so we can conclude
\textit{ii)}.
\end{pf}
\begin{prop}\label{induced-action}
Let $J:X\to S$ be a Hamiltonian space for the Manin pair
$(\mathfrak D, A)$. Then, there exists a Lie algebroid action of
$A\to S$ on $J:X\to S$.
\end{prop}
\begin{pf}
First, we will show that there exists a bundle map $\Phi :T^* X\to
A^*$ over $J:X\to S$. Given $\alpha \in T^*_xX$ one can define $\Phi
(\alpha ):A_{J(x)} \to \rr$ by
\[
\Phi (\alpha )(a )=\langle \alpha , u\rangle ,\mbox{ for }a\in
A_{J(x)},
\]
where, using \textit{i)} in Lemma \ref{vlemma}, $u\in T_xX$ is the
unique vector satisfying $(a, (u,0))\in F_{(J(x),x)}$.

Observe that, as the consequence of the definition of $\Phi :T^*X\to
A^*$, we have a map $\hat{}:\Gamma(A)\to \mathfrak X (X)$ such that
$(a(J(x)),(\hat{a}(x), 0))\in F_{(J(x),x)}$. Let us see that
$\hat{}:\Gamma(A)\to \mathfrak X (X)$, $a\mapsto \hat{a}$, indeed
defines a Lie algebroid action on $J:X\to S$.

First, since $F$ is compatible with the anchor ${\mathfrak D}\times
(TX\oplus T^*X)\to T(S\times X)$, $e=(d,(u,\alpha ))\mapsto
\tilde\rho (d)+u$, if we take $e=(a,(\hat{a},0))$, $a\in A_{J(x)}$,
we see that $\rho (a)+\hat{a}\in TQ$ if and only if
\[
J_*(\hat{a}(x))=\rho (a(J(x))).
\]
Moreover, from \textit{ii)} in Definition \ref{gds} in the
particular case when
$e_i|_{(J(x),x)}=((a_i(J(x)),0),(\hat{a}_i(x),0))$, $a_i\in \gm
(A)$, $x\in X$ and $i=1,2$, we get that
\[
\{ e_1, e_2\} = (\lcf a_1,a_2\rcf ,([\hat{a}_1,\hat{a}_2],0))
\]
and, as a consequence, $\{ e_1, e_2\}|_{(J(x),x)}\in F$ is
equivalent to
\[
\widehat{\lcf a_1, a_2\rcf} =[ \hat{a}_1 ,\hat{a}_2 ],\forall\,
a_1,a_2\in \gm (A).
\]
Therefore, we conclude our result.
\end{pf}
\begin{prop}\label{Dirac-X}
Let $J:X\to S$ be a Hamiltonian space for the Manin pair
$(\mathfrak D, A)$. Then,
\[
F(A)=\{ (u,\alpha )\in TX\oplus T^*X \, | \, (a,(u,\alpha ))\in F,
\, \exists\, a\in A\}
\]
is a Dirac structure on $TX\oplus T^*X$.
\end{prop}
\begin{pf}
First, let us see that $F(A)$ is an isotropic subbundle. If
$(u,\alpha ),(v,\beta )\in F(A)$ then there exists $a,b\in A$ such
that $ (a,(u,\alpha )), (b,(v,\beta ))\in F$. Using that $F\subset
\mathfrak D\oplus (TX\oplus T^*X)$ and $A\subset \mathfrak D$ are
isotropic, we get that
\[
0=\Big ( (a,(u,\alpha ))\Big |  (b,(v,\beta ))\Big )=(a\,|\,b)+\beta
(u)+\alpha (v)=\beta (u)+\alpha (v).
\]
Let $\alpha \in \mbox{Ker }\Phi _x\subset T_x^*X$. Then, from {\it
ii)} in Lemma \ref{vlemma}, we have that there exists $d\in
{\mathfrak D}_{J(x)}$ and $u\in T_xX$ such that $(d,(u,\alpha ))\in
F_{(J(x),x)}$. In addition, for any $a\in A_{J(x)}$, we have that
\[
0=\Big ( (d,(u,\alpha ))\Big |  (a,(\hat{a},0 ))\Big
)=(d\,|\,a)+\Phi (\alpha)(a)=(d\,|\,a).
\]
Since $A_{J(x)}$ is maximally isotropic, we deduce that $d\in
A_{J(x)}$ and, as a consequence, $(u,\alpha )\in F(A)_x$.

Suppose that $dim\mbox{ Ker }\Phi _x=i$. Take $\{ \alpha _1,\ldots
,\alpha _i\}$ a basis of $\mbox{Ker }\Phi _x$. Then, there exists
$u_1,\ldots, u_i$ such that $(u_1,\alpha _1),\ldots ,(u_i,\alpha
_i)\in F(A)_x$. On the other hand, take a basis $\{ v_1,\ldots
,v_{k-i} \}$ of $\mbox{ Im }(\Phi _x)^*$ ($k$ represents the
dimension of $X$). It is not difficult to prove that $\{ (u_1,\alpha
_1),\ldots ,(u_i,\alpha _i),(\hat{v}_1,0),\ldots ,(\hat{v}_{k-i},0)\}$ are
linearly independent in $F(A)_x$. Therefore $rank\, F(A)_x=dim\, X$

Finally, we prove that $\Gamma (F(A))$ is closed under the Courant
bracket $[\cdot ,\cdot ]_c$. Take $(X _1,\alpha _1),(X _2,\alpha
_2)\in \Gamma (F(A))$. There exist $a _1,a_2\in \Gamma (A)$ such
that $e_i|_{graph \,J}=(a _i,(X _i,\alpha _i))\in \Gamma (F)$. Since
$F$ is a generalized Dirac structure, we obtain that
\[
\gm (F)\ni \{ e_1,e_2\} |_{graph \,J}=([a_1,a_2],[(X_1,\alpha
_1),(X_2,\alpha _2)]_c ).
\]
Since $[a_1,a_2]\in \gm (A)$ ($A$ is a Dirac structure), we deduce
that $[(X_1,\alpha _1),(X_2,\alpha _2)]_c\in \gm (F(A))$.
\end{pf}
\subsection{Reduction}
In this Section we will show that, given a Hamiltonian space $J:X\to
S$, one can develop a reduction procedure in $X$, under some
regularity assumptions. More precisely,
\begin{them}
Let $J:X\to S$ be a Hamiltonian space for the Manin pair
$(\mathfrak D, A)$. Then, the orbit space $X/_{\widetilde{ }}$ is
a Poisson manifold
\end{them}
\begin{pf}
Let $J:X\to S$ be a Hamiltonian space for the Manin pair
$(\mathfrak D, A)$. Then, from Proposition (\ref{Dirac-X}), we
have that
\[
F(A)=\{ (u,\alpha )\in TX\oplus T^*X \, | \, (a,(u,\alpha ))\in F,
\, \exists\, a\in A\}
\]
is a Dirac structure on $TX\oplus T^*X$. Thus, according to
\cite{Co}, the set of admissible functions is endowed with a Poisson
algebra structure.

Next, let us describe the characteristic distribution of $F(A)$. If
$u_x\in F(A)\cap TX$ then there exists $a\in A_{J(x)}$ such that
\[
(a,(u_x,0 ))\in F_{(J(x),x)}
\]
and, as a consequence, $\hat{a}=u_x$. Therefore, the characteristic
distribution $F(A)\cap TX$ is just the distribution defined by the
action of the Lie algebroid $A$ over $J:X\to S$ and the manifold
$X/_{\widetilde{ }}\,$ is the reduction of the space by the action.

Since the set of admissible functions of $L$ is isomorphic to the
set of functions on the quotient $X/_{\widetilde{ }}\,$, we conclude
that $X/_{\widetilde{ }}\,$ is a Poisson manifold.
\end{pf}
\section{Hamiltonian spaces for Manin pairs and quasi-Lie bialgebroids}

\subsection{Hamiltonian spaces for quasi-Lie bialgebroids}
In this Section, we will recall the definition of a Hamiltonian
space for quasi-Lie bialgebroids and we will illustrate the notion
with a couple of examples.
\begin{defn}\label{Hamiltonian}
Let $(A,\delta ,\Omega )$ be a quasi-Lie bialgebroid over $S$. Then
$(X,\Pi _X)$ is a \textbf{Hamiltonian $A$-space} with
\textbf{momentum map $J:X\to S$} if there exists a Lie algebroid
action on $J:X\to S$ such that
\begin{equation}\label{quasi-action1}
[\Pi _X,J^\ast f]= \widehat{\delta f}, \, \forall\, f\in C^\infty
(S),
\end{equation}
\begin{equation}\label{quasi-action2}
[\Pi _X,\hat{a}]=\widehat{\delta a}, \ \forall\, a\in \Gamma (A) ,
\end{equation}
\begin{equation}\label{quasi-action3}
\half[\Pi _X ,\Pi _X]=\hat{\Omega}.
\end{equation}
\end{defn}

\begin{numex}
Let $(A,\delta ,\Omega )$ be a quasi-Lie bialgebroid over $S$. Then
$(S ,\Pi _S)$ is a Hamiltonian $A$-space with momentum map $Id:S\to
S$ where the Lie algebroid action is given by the anchor $\rho
:\Gamma (A)\to \mathfrak X(S)$ and $\Pi _S\in \mathfrak X^2(S)$ is
defined by
\[
\Pi _S (df_1,df _2)=- \langle \rho (\delta f_1),df_2\rangle ,\mbox{
for all }f_1,f_2\in C^\infty (S).
\]
We recall that Eq. (\ref{quasi-action3}) for this example has been
obtained in Proposition 4.8 in \cite{ILX}.
\end{numex}
\begin{numex}
Let $(A,A^*)$ be a Lie bialgebroid over $S$, that is, $(A,\delta
,\Omega )$ is a quasi-Lie bialgebroid with $\Omega =0$. Suppose that
$(X,\Pi _X)$ is a Hamiltonian $A$-space with momentum map $J:X\to
S$.

First, from Eq. (\ref{quasi-action3}), it is deduced that $(X,\Pi
_X)$ is a Poisson manifold. Moreover, it is not difficult to show
that Eqs. (\ref{quasi-action1}) and (\ref{quasi-action2}) are
equivalent to the equations defining a \textbf{Poisson action} of
the Lie bialgebroid $(A,A^*$) on the Poisson manifold $(X,\Pi _X)$
with momentum map $J:X\to S$ (see \cite{LWX}).

\end{numex}
\begin{numex}
Let $(\mathfrak g,F,\Omega )$ be a quasi-Lie bialgebra with
associated Manin quasi-triple $(\mathfrak d,\mathfrak g,\mathfrak
h)$. The Lie algebra action $\rho _X:\mathfrak g\to \mathfrak X(X)$
is a \textbf{quasi-Poisson action} of $(\mathfrak g,F,\Omega )$ on
$(X,\Pi_X)$ if
\[
\begin{array}{l}
[\rho _X(v),\Pi _X]=-\rho _X(F(v)), \mbox{ for }v\in \mathfrak
g,\\[7pt]
\displaystyle\half [\Pi _X,\Pi _X]=\rho _X(\Omega ),
\end{array}
\]
(see \cite{AK}). $(X,\Pi _X)$ is called a \textbf{quasi-Poisson
$\mathfrak g$-space}.
\end{numex}
\begin{numex}\label{Ex:quasi-Ham}
Let $(\mathfrak d,\mathfrak g,\mathfrak h)$ be a Manin quasi-triple
with corresponding quasi-Lie bialgebra $(\mathfrak g,F,\Omega )$.
Given $(X,\Pi_X)$ a quasi-Poisson $\mathfrak g$-space, an
equivariant map $J: X\to S$ is called a \textbf{momentum map} if
\begin{equation}\label{quasiPoissonAction3}
\Pi _X^\sharp (J^\ast \theta _s)=-(\lambda (\theta _s))_X,\mbox{ for
}\theta _s\in T^\ast _sS,
\end{equation}
where $\mathfrak g$ acts on $S$ by the dressing action, and $\lambda
:T^\ast _sS\to \mathfrak g$ is the dual map of the infinitesimal
dressing action $\mathfrak g^\ast \cong \mathfrak h \to T_sS$, i.e.,
\[
\langle \lambda (\theta _s),\eta \rangle = \langle \theta _s,\eta
_S(s) \rangle, \forall\, \theta _s \in T^\ast _sS\mbox{ and }\eta
\in \mathfrak h.
\]
In this case, $X$ is called a \textbf{quasi-hamiltonian space} (see
\cite{AK}).
\end{numex}

\subsection{Hamiltonian spaces for Manin pairs and quasi-Lie
bialgebroids} \label{split-case} Now, we will study Hamiltonian
spaces for Manin pairs $(\mathfrak D, A)$ in the particular case
when $A$ admits an isotropic complement in $\mathfrak D$.

\medskip
Let $J:X\to S$ be a Hamiltonian space for the Manin pair
$(\mathfrak D, A)$. If we choose an isotropic complement $B$ of
$A$ in ${\mathfrak D}$, then we can identify $B\cong A^*$ and
${\mathfrak D}=A\oplus B\cong A\oplus A^*$.

\medskip
First, since $(a, (\hat{a} , 0))\in F_{(J(x),x)}$, for any $a\in
A_{J(x)}$, the isotropy property of $F$ implies that any $e\in F$
must be of the form
\begin{equation}\label{decomposition}
e\in F \Longrightarrow e= ((a,-\Phi (\alpha )),(u,\alpha )),\mbox{
for }a\in A, \alpha \in T^*X\mbox{ and }u\in TX,
\end{equation}
where $\Phi$ is the dual of the action given by Proposition
\ref{induced-action}.

\medskip
Next, we will define a bivector $\Pi _X\in \mathfrak X^2(X)$. Given
$\alpha \in T^*X$, we characterize $\Pi ^\sharp _X (\alpha )\in TX$
by
\[
\beta (\Pi ^\sharp _X (\alpha ))=\Phi (\beta )(a)-\beta
(u),\forall\, \beta \in T^*X ,
\]
where, using Lemma \ref{vlemma} \textit{ii)} and Eq.
(\ref{decomposition}), $u\in TX$ and $a\in A$ satisfy $((a,-\Phi
(\alpha )),(u,\alpha ))\in F$. Let us see that $\Pi ^\sharp _X$ is
well defined. If $((a',-\Phi (\alpha )),$ $(u',\alpha ))\in F$ then
$((a'-a,0),(u'-u,0 ))\in F$. From {\it ii)} in Definition
\ref{Hamiltonian-space}, we see that $u'-u=\widehat{a'-a}$. As a
consequence, $\Phi (\beta )(a)-\beta (u) =\Phi (\beta )(a')-\beta
(u')$.

On the other hand, from the isotropy property of $F$ we have that,
for any $\alpha ,\beta \in T^*X$,
\[
0= \Big ( ((a,-\Phi (\alpha )),(u,\alpha )) \Big | ((b,-\Phi (\beta
)),(v,\beta ))  \Big ) = \beta (\Pi ^\sharp _X(\alpha )) +\alpha
(\Pi ^\sharp _X (\beta )).
\]
That is, $\Pi ^\sharp _X$ is skew-symmetric.

\medskip
Summing up, if we choose an isotropic complement $B$ of $A$ in
${\mathfrak D}$ then a Hamiltonian space $J:X\to S$ for the Manin
pair $({\mathfrak D},A)$ can be described by a Lie algebroid
action $\hat{ }:\Gamma (A)\to \mathfrak X(X)$ and a bivector $\Pi
_X$ such that the generalized Dirac structure is given by
\[
F\kern-1pt=\kern-1pt\{ ((a,-\Phi (\alpha )),(\hat{a}-\Pi _X^\sharp
(\alpha ),\alpha )) \,|\, \forall\, a\in A_{J(x)}, \,\alpha \in
T_x^*X\} \kern-1pt\subset\kern-1pt (A\oplus A^*)\times (TX\oplus
T^*X).
\]
Moreover,
\begin{them}\label{action-qlb}
Let $(A,\delta ,\Omega )$ be a quasi-Lie bialgebroid on $S$ and
$(X,\Pi _X)$ such that there exists a bundle map $\Phi :T^\ast
X\to A^*$ over $J:X\to S$. Then, $X$ is a Hamiltonian $A$-space
with momentum map $J:X\to S$ if and only if $J:X\to S$ is a
Hamiltonian space for the Manin pair $(A\oplus A^*,A)$, where the
generalized Dirac structure $F\to Q=graph\, J$ is given by
\[
F=\{ ((a,-\Phi (\alpha )),(\hat{a}-\Pi _X^\sharp (\alpha ),\alpha ))
\,|\, \forall\, a\in A_{J(x)}\mbox{ and }\alpha \in T_x^*X\} .
\]
\end{them}
\begin{pf}
First, we will check that $F$ is compatible with the anchor
$\tilde{\rho }:(A\oplus A^\ast)\times (TX\oplus T^*X)\to T(S\times
X)$, $((a,\xi ),(u,\alpha ))\mapsto \rho (a)+\rho _*(\xi )+u$.

\medskip \textit{First case}. $e=((a,0),(\hat{a},0))$, $a\in
A_{J(x)}$.

It is trivial to see that $\tilde{\rho} (e)=\rho (a)+\hat{a}\in TQ$
if and only if
\begin{equation}\label{1}
J_*(\hat{a}(x))=\rho (a(J(x))).
\end{equation}

\medskip
\textit{Second case}. $e=((0,-\Phi (\alpha )),(-\Pi ^\sharp_X(\alpha
) ,\alpha ))$, $\alpha \in T^*_xX$.

In this case, $\tilde{\rho} (e)\in TQ$ if and only if
\begin{equation}\label{2'}
\rho _*(\Phi (\alpha ))=J_* (\Pi ^\sharp_X(\alpha )).
\end{equation}
Consider $g\in C^\infty (M)$. Then,
\[
\begin{array}{l}
J_* (\Pi ^\sharp_X(\alpha ))(g)=\langle [\Pi _X,J^*g],\alpha
\rangle,\mbox{ and }\\[5pt]
\rho _*(\Phi (\alpha ))(g)=\langle \widehat{\delta g}(x),\alpha
\rangle .
\end{array}
\]
Thus, Eq. (\ref{2'}) is equivalent to
\begin{equation}\label{2}
\widehat{\delta g}=[\Pi _X,J^*g], \,\forall\, g\in C^\infty (M).
\end{equation}
Now, let us check condition \textit{ii)} in Definition \ref{gds}.
We will distinguish three different cases here.

\medskip
\textit{First case}.
$e_i|_{(J(x),x)}=((a_i(J(x)),0),(\hat{a}_i(x),0))$, $a_i\in \gm
(A)$, $x\in X$ and $i=1,2$.
\[
\{ e_1, e_2\} = ((\lcf a_1,a_2\rcf ,0),([\hat{a}_1,\hat{a}_2],0))
\]
Therefore, $\{ e_1, e_2\}|_{(J(x),x)}\in F$ is equivalent to
\begin{equation}\label{3}
\widehat{\lcf a_1, a_2\rcf} =[ \hat{a}_1 ,\hat{a}_2 ],\forall\,
a_1,a_2\in \gm (A).
\end{equation}

\medskip
\textit{Second case}. $e_i|_{(J(x),x)}=((0, -\Phi (\alpha _i(x))
),(-\Pi ^\sharp _X (\alpha _i(x)),\alpha _i(x) ))), \forall
x\in X$ and $i=1,2$. Here $\alpha _i\in
\Omega ^1(X)$ is chosen so that $\Phi (\alpha _i)$
is a well defined section in $\Gamma (A^*)$.

\[
\{ e_1, e_2\} =  (((\Phi (\alpha _1)\wedge \Phi (\alpha _2))\per
\Omega , \lcf \Phi (\alpha _1),\Phi (\alpha _2)\rcf _* ),([\Pi
^\sharp _X(\alpha _1) , \Pi ^\sharp _X(\alpha _2)], -\lcf \alpha
_1,\alpha _2\rcf _{\Pi _X} )),
\]
where $\lcf \cdot ,\cdot \rcf _{\Pi _X}$ is the bracket on $\Omega
^1(X)$ associated with the bivector field $\Pi _X$. $\{ e_1,
e_2\}|_Q\in \gm (F)$ if and only if
\[
\begin{array}{l}
a)\; \Phi (\lcf \alpha _1,\alpha _2 \rcf _{\Pi _X} )= \lcf \Phi
(\alpha _1),\Phi (\alpha _2)\rcf _* ,\\[8pt]
b)\; [\Pi ^\sharp _X(\alpha _1) , \Pi ^\sharp _X(\alpha _2)] - \Pi
^\sharp _X (\lcf \alpha _1,\alpha _2\rcf _{\Pi _X} )= \Big ( (\Phi
(\alpha _1)\wedge \Phi (\alpha _2))\per \Omega \Big )\widehat{
}\,.
\end{array}
\]
Using the identity $ [\Pi ^\sharp _X(\alpha _1) , \Pi ^\sharp
_X(\alpha _2)] - \Pi ^\sharp _X (\lcf \alpha _1,\alpha _2\rcf _{\Pi
_X} )= (\alpha _1\wedge\alpha _2) \per \half [\Pi _X,\Pi _X]$, we
see that $b)$ is equivalent to
\begin{equation}\label{4}
\half [\Pi _X,\Pi _X] =\hat{\Omega}.
\end{equation}
On the other hand, from Eqs. (\ref{bracket+anchor}) and
(\ref{2'}), we observe that $a)$ is equivalent to
\begin{equation}\label{5}
\widehat{\delta a}=[\Pi _X,\hat{a}], \forall\, a\in \gm (A).
\end{equation}

\medskip
\textit{Third case}. $e_1|_{(J(x),x)}=((a(J(x)),0),(\hat{a}(x),0))$,
$a\in \gm (A)$, $x\in X$ and $e_2|_{(J(x),x)}=((0, -\Phi (\alpha
(x)) ),(-\Pi ^\sharp _X (\alpha (x)),\alpha (x)))$, $\alpha \in
\Omega ^1(X)$, $x\in X$.

If we compute the bracket $\{ e_1,e_2\}$, we obtain
\[
\{ e_1,e_2\} = (( i_{\Phi (\alpha )}\delta a , -{\cal L}_a\Phi
(\alpha )),(-[\hat{a},\Pi ^\sharp _X (\alpha )],{\cal
L}_{\hat{a}}\alpha )).
\]
Thus, $\{ e_1,e_2\}|_{(J(x),x)}\in F$ for any $x\in X$ if and only
if
\[
\begin{array}{l}
c)\; {\cal L}_a(\Phi (\alpha ) )= \Phi ( {\cal L}_{\hat{a}}
\alpha ), \\[8pt] d)\; [\hat{a},\Pi ^\sharp _X (\alpha )]= -\widehat{i_{\Phi
(\alpha )}\delta a}+ \Pi ^\sharp _X( {\cal L}_{\hat{a}} \alpha ).
\end{array}
\]
Condition $c)$ is immediately deduced from the definition of the
Lie derivative and Eqs. (\ref{1}) and (\ref{3}).

On the other hand, for any $\Pi \in \mathfrak X^2 (X)$, $X\in
\mathfrak X(X)$ and $\alpha _1,\alpha _2\in \Omega ^1(X)$ we know
that
\[
\begin{array}{rcl}
[X,\Pi ](\alpha _1 ,\alpha _2 )&=& X(\Pi (\alpha _1 ,\alpha _2
))-\langle \alpha _2 , \Pi ^\sharp ({\cal L}_X\alpha _1 ) \rangle
-\langle {\cal L}_X\alpha _2,
\Pi ^\sharp (\alpha _1 )\rangle \\
&=& -\langle \alpha _2 ,\Pi ^\sharp ({\cal L}_X\alpha _1 ) \rangle
+\langle \alpha _2 ,[X,\Pi ^\sharp (\alpha _1 ) ]\rangle.
\end{array}
\]
Since $\widehat{i_{\Phi (\alpha )}\delta a}=i_\alpha
\widehat{\delta a}$, we have that $d)$ is equivalent to Eq.
(\ref{5}).

Summing up, we have that Eqs. (\ref{1}) and (\ref{3}) are equivalent
to a Lie algebroid $A$-action on $J:X\to S$ and Eqs. (\ref{2}),
(\ref{4}) and (\ref{5}) are just the conditions in Definition
\ref{Hamiltonian}.
\end{pf}
\begin{rmk}
If we choose an isotropic complement $B$ of $A$, then there exists a
bivector field $\Pi _X$ and
\[
F(A)= \{ (\hat{a}-\Pi ^\sharp_X(\alpha ),\alpha )\, | \, \alpha \in
T^*X ,\, a\in A,\, \Phi (\alpha )=0 \}.
\]
A direct computation shows that the Poisson bivector on the the
orbit space is just the projection of $-\Pi _X$.
\end{rmk}
\begin{numex}
Let $(\mathfrak g, F,\Omega )$ be a Lie quasi-bialgebra and
$(\mathfrak d=\mathfrak g\oplus \mathfrak g^*,\mathfrak g, \mathfrak
g^*)$ be the corresponding Manin quasi-triple (see Eq.
(\ref{d-bracket}) in Example \ref{Manin-pair}). Using Theorem
\ref{action-qlb} we recover the result obtained in \cite{BCS}, that
is, $(X,\Pi _X)$ is a quasi-Poisson $\mathfrak g$-space if and only
if
\[
\{ ((\rho_X(u)-\Pi _X^\sharp (\alpha ),\alpha ),(u,-\Phi (\alpha )))
\,|\, \forall\, u\in \mathfrak g, \,\alpha \in T_x^*X\}
\]
is a Dirac structure in $\mathfrak d \times (TX\oplus T^*X)$. In
addition, the characteristic distribution of $F(\mathfrak g\times
S)$ is just the image of the action of $\mathfrak g$ on $X$ and the
Poisson structure on the reduced space is just given by the
restriction of $\Pi _X$. This was first proved in \cite{AK} (see
Theorem 4.2.2).
\end{numex}

\begin{numex}\label{Ex:quasi-ham}
Let $(X,\Pi _X)$ be quasi-Hamiltonian space for a Manin
quasi-triple $(\mathfrak d,\mathfrak g,\mathfrak h)$ with momentum
map $J:X\to S$. From the Manin quasi-triple $(\mathfrak
d,\mathfrak g,\mathfrak h)$ one can construct the Manin pair
$(\mathfrak d\times S, \mathfrak g\times S)$ which admit an
splitting $\mathfrak d\times S\cong (\mathfrak g\times S)\oplus
(\mathfrak h\times S)$. Using Example \ref{Ex:quasi-ham} and
Theorem \ref{action-qlb}, we have that a quasi-hamiltonian
$\mathfrak g$-space $(X,-\Pi _X)$ with momentum $J:X\to S$ is
equivalent to a Hamiltonian space $J:X\to S$ for the Manin pair
$(\mathfrak d\times S,\mathfrak g\times S)$. In this case, the
Dirac structure $F(\mathfrak g \times S)$ is given by
\[
F(\mathfrak g\times S)= \{ (\rho _X(v)+\Pi ^\sharp_X(\alpha ),\alpha
)\, | \, \alpha \in T^*X ,\, v \in \mathfrak g,\,  \Phi (\alpha )=0
\}.
\]
\end{numex}

\begin{numex}
Let $(\mathfrak d,\mathfrak g,\mathfrak h)$ be a Manin quasi-triple
with corresponding quasi-Lie bialgebra $(\mathfrak g,F,\Omega )$.
Moreover, suppose that $(X,\Pi_X)$ is a quasi-Hamiltonian space with
momentum map $J: X\to S$.

If we consider the quasi-Lie bialgebroid $(\mathfrak g\times
S,\delta ,\Omega )$ (see Example \ref{Manin-pair}), it is easy to
show that $(X,\Pi _X)$ is a quasi-Hamiltonian space for $(\mathfrak
g,F,\Omega)$ if and only if $(X,-\Pi _X)$ is a Hamiltonian
$(\mathfrak g\times S)$-space for the quasi-Lie bialgebroid
$(\mathfrak g\times S,\delta ,\Omega )$. As a consequence, $(X,\Pi
_X)$ is a quasi-Hamiltonian space for $(\mathfrak g,F,\Omega)$ with
momentum map $J:X\to S$ if and only if $J:X\to S$ is a  Hamiltonian
space $J:X\to S$ for the Manin pair $(\mathfrak d\times S,\mathfrak
g\times S)$. Moreover, an straightforward computation shows that the
characteristic distribution of $F(\mathfrak g\times S)$ is just the
image of the action of $\mathfrak g$ on $X$ and the Poisson
structure on the reduced space is just given by $\Pi _X$. This was
first proved in \cite{AK} (see Theorem 4.2.2).
\end{numex}
\begin{numex}
If the isotropic complement $H$ is a Dirac structure, we know that
$(A, A^*)$ is a Lie bialgebroid and, thus, $(X,\Pi _X)$ is a Poisson
manifold . Using the Courant algebroid isomorphism $TX\oplus _{\Pi
_X}T^*X\to TX\oplus T^*X$, $(v,\alpha)\mapsto (\Pi ^\sharp _X
(\alpha )+v,\alpha )$ we get that
\[
F=\{ ((a,-\Phi (\alpha )),(\hat{a},\alpha )) \,|\, \forall\, a\in
A_{J(x)}\, \mbox{and} \,\alpha \in T_x^*X\} \subset (A\oplus
A^*)\times (TX\oplus T^* X)
\]
is a generalized Dirac structure in $(A\oplus A^*)\times (TX\oplus
_{\Pi _X}T^*X)$. As a consequence (see Example \ref{ej-gds}), we
conclude that the dual of the action $\Phi :T^\ast X\to A^\ast$ is a
Lie bialgebroid morphism between $(T^*X,TX)$ and $(A^*,A)$.
\end{numex}
A particular situation of the previous example is the following one.
\begin{numex}
Let $J:X\to S$ be a Poisson map. Then, using that it induces a Lie
bialgebroid morphism, it is a Hamiltonian space for the Manin pair
$(T^*S\oplus _{\Pi _S} TS, T^\ast S)$. In this case, the Dirac
structure $F(T^*S)$ is given by
\[
F(T^*S)= \{ (\Pi ^\sharp_X(J^\ast \beta -\alpha),\alpha )\, | \,
\alpha \in T^*X ,\, \beta \in T^*S,\,  J_\ast \Pi _X^\sharp (\alpha
)=0 \}\subset TX\oplus T^*X,
\]
the characteristic distribution is
\[
{\mathfrak l}=F(T^*S)\cap TX= \{ \Pi ^\sharp _X (J^\ast \beta ) \,|
\, \beta \in T^\ast X \}.
\]
and the Poisson structure on $X/{\mathfrak l}$ is just the
restriction of $\Pi _X$.
\end{numex}

\bigskip
{\bf Acknowledgments.} Iglesias Ponte's research was partially
supported by grants MTM2007-62478 and S-0505/ESP/0158 of the CAM,
and Xu's research was partially supported by NSF grants
DMS-0605725 and DMS-0801129 \&  NSA grant H98230-06-1-0047.
Iglesias Ponte thanks CSIC for a ``JAE-Doc" research contract and
Mathematics Department at Penn State University, where part of
this project was done. We also thank Zhou Chen for many useful
comments.

\end{document}